\documentclass[psfig,a4paper,12pt,reqno]{amsart}
\usepackage{amsfonts,hyperref}
\usepackage{amssymb,amsmath,amsbsy,array,latexsym,color}

 \theoremstyle{plain}
\begingroup
\newtheorem{theorem}    {Theorem}                       [section]
\newtheorem{lemma}      [theorem]        {Lemma}
\newtheorem{proposition}[theorem]        {Proposition}
\newtheorem{corollary}  [theorem]        {Corollary}
\newtheorem{example}    [theorem]        {Example}
\theoremstyle{definition}
\newtheorem{remark}     [theorem]        {Remark}
\newtheorem{definition} [theorem]        {Definition}

\numberwithin{equation}                                 {section}

\title{On sublaplacians of sub-Riemannian manifolds}
\author{Kang-Hai Tan}
\address{Institute of Mathematics\\Fudan University, 200433\\ Shanghai, PRC}

\email{tankanghai2000@yahoo.com.cn,khtan@fudan.edu.cn}

\begin{document}

\begin{abstract}
 \small{In this note we address a notion of sublaplacians
 of sub-Riemannian manifolds. In particular for fat sub-Riemannian manifolds
 we answered the sublaplacian question proposed by R. Montgomery in
 \cite[p.142]{Montgomery}
}
\end{abstract}

\maketitle

\section{Introduction}
In the literature the sub-Riemannian analogue of Riemannian
Laplacians is H\"{o}rmander sum of squares of vector fields, see
\cite{Hor,RS,Bi,VST,Str,GV1,Taylor1,Taylor2,Ge,Ru1,Ru2} and
references therein. Let $M$ be a smooth manifold of dimension $m$
endowed with a smooth distribution (horizontal bundle) $\Sigma$ of
dimension $k$ with $k<m$. If we a prior equip $\Sigma$ with an inner
product $g_{c}$ (sub-Riemannian metric), we call $(M,\Sigma,g_{c})$
a sub-Riemannian manifold with the sub-Riemannian structure
$(\Sigma,g_{c})$. If $\Sigma$ is integrable, it is just the
Riemannian geometry. We will assume $\Sigma$ is  not integrable. A
piecewise smooth curve $\gamma(t),t\in{[a,b]}$ in $M$ is horizontal
if $\dot{\gamma}(t)\in{\Sigma_{\gamma(t)}}$ a.e. $t\in{[a,b]}$. The
\textit{length} $\ell(\gamma)$ of the horizontal curve
$\gamma(t),t\in{[a,b]}$ is the integral
$\int_{a}^{b}g_{c}(\dot{\gamma}(t),\dot{\gamma}(t))dt$. Denote by
$\Sigma_{i}$ the set of all vector fields spanned by all commutators
of order $\leq i$ of vector fields in $\Sigma$  and let
$\Sigma_{i}(p)$ be the subspace of evaluations  at $p$ of all vector
fields in $\Sigma_{i}$.  We call $\Sigma$ satisfies the Chow or
H\"{o}rmander condition if for any $p\in{M}$, there exists an
integer $r(p)$ such that $\Sigma_{r(p)}(p)=T_{p}M$ (the least such
$r$ is called the degree of $\Sigma$ at $p$). If moreover
$\Sigma_{i}$ is of constant dimension for all $i\leq r$, $\Sigma$
and also $(M,\Sigma,g_c)$ are called regular. If $M$ is connected
and $\Sigma$ satisfies the H\"{o}rmander condition, the Chow
connectivity theorem asserts that there exists at least one
piecewise smooth horizontal curve connecting two given points (see
\cite{Chow,Sub-Rie,Montgomery}), and thus $(\Sigma,g_{c})$ yields a
metric (called \textit{Carnot-Carath\'eodory distance}) $d_{cc}$ by
letting $d_{cc}(p,q)$ as the infimum among the lengths of all
horizontal curves joining $p$ to $q$. Let $\{X_{1},\cdots,X_{k}\}$
be an orthonormal basis of $\Sigma$.  The H\"{o}rmander operator is
$\square=\sum_{i=1}^kX_i^2+X_0$ where $X_0\in{\Gamma(\Sigma)}$ is a
horizontal vector field. It is easy to see that  the operator
$\square$ in general depends on the choice of orthonormal bases.
Thus $\square$ is not intrinsic to the sub-Riemannian structure
$(\Sigma,g_c)$. Recall that the Riemannian Laplacian on a Riemannian
manifold is a Riemannian invariant. Montgomery in
\cite[p.142]{Montgomery} proposed the question whether there exists
a canonical sublaplacian in the sub-Riemannian case. As observed by
\cite{Montgomery}, this question is equivalent to the existence of a
canonical measure $\mu$: the canonical sublaplacian
$\widehat{\Delta}$ and $\mu$ should satisfy
$$
-\int_M(\widehat{\Delta}f)gd\mu=\int_Mg_c(\nabla^{\mathcal{H}}f,\nabla^{\mathcal{H}}g)d\mu
$$
for any smooth functions $f,g$ with compact support, where
$\nabla^{\mathcal{H}}f$ is the horizontal gradient of $f$:
$g_c(\nabla^{\mathcal{H}}f,X)=Xf$ for any horizontal vector field
$X\in{\Gamma(\Sigma)}$.

The importance of the study of sublaplacians on sub-Riemannian
manifolds lies in the conjectured close relationship between
spectral asymptotics of sublaplacians and sub-Riemannian geodesics
(in particular singular curves), see \cite{Montgomery,Montgomery1}
and references therein. This note is devoted to a rudimental study
of sublaplacians. We will give a notion of \textbf{sublaplacians}
for general sub-Riemannian manifolds, or \textbf{the sublaplacian}
for \textbf{fat} sub-Riemannian manifolds. We recall that a
sub-Riemannian manifold $(M,\Sigma,g_c)$ is \textbf{fat} if $\Sigma$
is strong-bracket generating, that is, for each $p\in{M}$ and each
nonzero horizontal vector $v\in{\Sigma_p}$ we have
$$
\Sigma_p+[V,\Sigma]_p=T_pM
$$
where $V$ is any horizontal extension of $v$. This notion is defined
using the truncated connection (called horizontal connection). Given
a complement $\Sigma^{\prime}$ of $\Sigma$:
$TM=\Sigma\bigoplus\Sigma^{\prime}$, then with respect to this
decomposition there exists a unique horizontal connection $D$ on
$\Sigma$ such that for any $X,Y,Z\in{\Gamma(M)}$
\begin{enumerate}
\item $Xg_c(Y,Z)=g_c(D_XY,Z)+g_c(Y,D_XZ)$
\item $D_XY-D_YX=[X,Y]^{\mathcal{H}},$
\end{enumerate}
where $[X,Y]^{\mathcal{H}}$ denotes the projection on $\Sigma$ of
$[X,Y]$, see \cite{TY,Ge} for details. We define
$$
\Delta^{\mathcal{H}}:=\mathrm{div}^{\mathcal{H}}\circ\nabla^{\mathcal{H}}
$$
where the horizontal divergence operator
 is defined as
 $\mathrm{div}^{\mathcal{H}}X:=\sum_{i=1}^kg_c(D_{X_i}X,X_i)$ for
 $X\in{\Gamma(\Sigma)}$ and an orthonormal basis $\{X_i\}_{i=1}^k$.
 The operator $\Delta^{\mathcal{H}}$ depends only on $(\Sigma,g_c)$
 and the decomposition. In many cases such as nilpotent groups, contact manifolds, principal bundles with
 connections, and Riemannian submersions,
 there exists a `natural'
 decomposition of $TM$. In particular when the sub-Riemannian metric
 $g_c$ is the projection on $\Sigma$ of a Riemannian metric $g$,
 $TM$ can be orthogonally decomposed as
 $TM=\Sigma\bigoplus\Sigma^{\prime}$. Conversely, given a
 decomposition, we always can extend $g_c$ to a Riemannian metric
 $g$ such that the decomposition is orthogonal. Note that such
extensions are not unique. Let $g$ be any such
 extension. Then for any $X,Y\in{\Gamma(\Sigma)}$
 $$
 D_XY=\mathcal{P}(\nabla_XY)
 $$
where $\nabla$ is the Riemannian connection of $g$ and $\mathcal{P}$
denotes the projection on $\Sigma$, see \cite{TY}.

As stated above the defined $\Delta^\mathcal{H}$ depends on the
splitting of the tangent bundle. This makes the problem delicate. We
remark that for some cases such as nilpotent groups with grading Lie
algebra  and contact Riemannian manifolds, there is a canonical
notion of the sublaplacian.  Recall that the sub-Riemannian geometry
of $(M,\Sigma,g_c)$, i.e., the geometry of $(M,d_{cc})$, depends
only on the sub-Riemannian structure $(\Sigma,g_c)$, not on
complements of $\Sigma$ or extensions of $g_c$. Our motivation to
study sublaplacians or weakly convex function on sub-Riemannian
manifolds is to extract information about sub-Riemannian geometry as
much as possible by exploring functions or invariants defined
intrinsicly. What should a canonical complement of $\Sigma$, and
then a canonical orthogonal extension of $g_c$ be? For instance, in
our opinion it is desirable (under some conditions imposed on
$\Sigma$ and topology of $M$) to find a complement of (regular)
$\Sigma$ and then to select an extension $g$ of $g_c$ such that the
Riemannian measure of $g$ is just the Hausdorff $Q-$measure of
$d_{cc}$, where $Q=\sum_{i=1}^ri(\dim(\Sigma_i)-\dim(\Sigma_{i-1}))
( \Sigma_{0}=\emptyset,\Sigma_{1}=\Sigma$) is the Hausdorff
dimension of $(M,d_{cc})$. The following statement is our starting
point.

\begin{theorem}\label{main2}Let $(M,\Sigma,g_c)$ be a sub-Riemannian manifold.
Then there exists a complement $\Sigma^\prime$ of $\Sigma$,
$TM=\Sigma\bigoplus\Sigma^{\prime}$, such that for this
decomposition there is an orthogonal extension $g$ of $g_c$ and an
orthonormal basis $\{T_1^\prime,\cdots,T^\prime_{m-k}\}$ of
$\Sigma^\prime$ satisfying
\begin{equation}\label{curvature}
\mathcal{P}(\nabla_{T^\prime_\beta}T^\prime_\beta)=0,\quad\beta=1,\cdots,m-k
\end{equation}
where $\nabla$ is the Riemannian connection of $g$. Moreover, if for
any local frame of $TM$, $\{X_1,\cdots,X_k,$ $T_1,\cdots,T_{m-k}\}$,
where $\{X_1,\cdots,X_k\}$ is a basis of $\Gamma(\Sigma)$, the
matrix $[C_{ij}^\beta]$ is invertible for $\beta=1,\cdots,m-k,$
where $C_{ij}^\beta$ are the coefficients satisfying
\begin{equation}\label{coefficients}
[X_i,X_j]=\sum_{a=1}^kC_{ij}^a X_a+\sum_{\beta=1}^{m-k}C_{ij}^\beta
T_\beta,
\end{equation}

then such complement is unique.
\end{theorem}
In general, the condition in Theorem \ref{main2} guaranteeing the
uniqueness of the complement is very strong. For most cases it is
impossible. In fact this condition is just the strong-bracket
generating condition for $\Sigma$, see Proposition
\ref{strongbracket}. We recall that contact structures of  contact
manifolds are fat.

The following theorem motivates our definition of sublaplacians, see
Definition \ref{defn}.
\begin{theorem}\label{main1}
Let $(M,\Sigma,g_c)$ be a sub-Riemannian manifold. Assumme
$TM=\Sigma\bigoplus\Sigma^{\prime}$ be a splitting and $g$ an
orthogonal extension of $g_c$. Denote by $d\mathrm{vol}$ the
Riemannian measure of $g$ and by $\mathrm{H}^\bot$ the mean
curvature of $\Sigma^{\prime}$. Let $u$ be a positive, smooth
function on $M$. Set $d\mu=ud\mathrm{vol}$. Then
$$
-\int_M(\Delta^{\mathcal{H}}e)fd\mu=\int_Mg_c(\nabla^{\mathcal{H}}e,\nabla^{\mathcal{H}}f)d\mu
$$
holds for any $e,f\in{C_0^{\infty}(M)}$ if and only if $\ln u$ is
the horizontal potential of $\mathrm{H}^\bot$, i.e.,
\begin{equation}\label{characterization}
\nabla^{\mathcal{H}}(\ln u)=\mathrm{H}^\bot.
\end{equation}
\end{theorem}
Here the mean curvature $\mathrm{H}^\bot$ of $\Sigma^{\prime}$ is
defined as
$$
\mathrm{H}^\bot:=\sum_{\beta=1}^{m-k}\mathcal{P}
\left(\nabla_{T_\beta}T_\beta\right)=\sum_{\beta=1}^{m-k}\sum_{i=1}^k
g\left(\nabla_{T_\beta}T_\beta,X_i\right)X_i,
$$
where $\{X_i\}_{i=1}^k,\{T_\beta\}_{\beta=1}^{m-k}$ are orthonormal
bases of $\Sigma$, $\Sigma^{\prime}$ respectively. We note that for
given $\Sigma^\prime$ and $g$, equation \eqref{characterization} is
not soluble in general. But by Theorem \ref{main2}, we always can
choose an orthogonal extension of $g_c$ such that
\eqref{characterization} is soluble for $\mathrm{H}^\bot=0$, if
$(M,\Sigma,g_c)$ satisfies the assumption in Theorem \ref{main2}.

The paper is organized as follows. In the next section after proving
Theorem \ref{main2}, \ref{main1}, we give the definition of
sublaplacians. Several canonical examples are given. It turns out
that our definition is compatible with the canonical sublaplcians in
the literature. At the end of Section \ref{sec2}, a Hopf type
theorem is proven for closed sub-Riemannian manifolds. The last
section is devoted to the closed eigenvalue problem on compact
sub-Riemannian manifolds.
\section{Proofs, examples, and basic properties of
sublaplacians}\label{sec2}
\begin{proof}[\textbf{Proof of Theorem \ref{main2}}]Let $\{X_1,\cdots,X_k,T_1,\cdots,T_{m-k}\}$
 be any local frame of $TM$,
where $\{X_1,\cdots,X_k\}$ is an orthonormal basis of
$\Gamma(\Sigma)$. Extend $g_c$ to a Riemannian metric $\bar{g}$ such
that $\{T_\beta\}_{\beta=1}^{m-k}$ is orthonormal. Denote by
$\bar{\nabla}$ the Riemannian connection of $\bar{g}$. Then for
$\beta=1,\cdots,m-k,$
$$
\bar{F}_\beta:=\mathcal{P}(\bar{\nabla}_{T_\beta}T_\beta)=\sum_{i=1}^k\bar{g}(T_\beta,[X_i,T_\beta])X_i.
$$
Let $\Sigma^\prime$ be any complement of $\Sigma$ and
$\{\bar{T}_\beta\}_{\beta=1}^{m-k}$ be a local basis of
$\Sigma^\prime$. Then
$$
 \bar{T}_\beta=\sum_{i=1}^k\bar{A}^i_\beta
 X_i+\sum_{\alpha=1}^{m-k}B^\alpha_\beta T_\alpha
$$
for smooth functions $\bar{A}^i_\beta$ and $B^\alpha_\beta$. Since
$\{\bar{T}_\beta\}_{\beta=1}^{m-k}$ is a basis of $\Sigma^\prime$,
the matrix $[B_\beta^\alpha]$ is invertible. Denote by
$[K_\beta^\alpha]$ the inverse matrix of $[B_\beta^\alpha]$. Setting
\begin{align*}
T^\prime_\beta:&=\sum_{\alpha=1}^{m-k}K_\beta^\alpha\bar{T}_\alpha\\
&=\sum_{\alpha=1}^{m-k}\sum_{i=1}^{k}K_\beta^\alpha\bar{A}^i_\alpha
X_i+\sum_{\alpha=1}^{m-k}\sum_{\gamma=1}^{m-k}K_\beta^\alpha
B^\gamma_\alpha T_\gamma,
\end{align*}
i.e.,
$$
T^\prime_\beta=\sum_{i=1}^kA^i_\beta X_i+
 T_\beta, \quad\mathrm{where}\quad
 A^i_\beta=\sum_{\alpha=1}^{m-k}K_\beta^\alpha\bar{A}^i_\alpha,
$$
then $\{T^\prime_\beta\}_{\beta=1}^{m-k}$ is also a basis of
$\Sigma^\prime$. Now extend $g_c$ to a Riemannian metric $g$ such
that $\{T^\prime_1,\cdots,T^\prime_{m-k}\}$ is orthonormal with
respect to $g$. Noting that
\begin{align*}
[X_i,T_\beta^\prime]&=[X_i,\sum_{j=1}^kA^j_\beta
X_j+T_\beta]\\&=\sum_{j=1}^kA^j_\beta[X_i,X_j]+(X_iA^j_\beta)X_j+[X_i,T_\beta],
\end{align*}
we have for $\beta=1,\cdots,m-k$
\begin{align*}
\mathcal{P}(\nabla_{T^\prime_\beta}T^\prime_\beta)&=\sum_{i=1}^kg(T^\prime_\beta,[X_i,T^\prime_\beta])X_i\\
&=\sum_{i=1}^kg\left(T^\prime_\beta,\sum_{j=1}^kA^j_\beta[X_i,X_j]+(X_iA^j_\beta)X_j+[X_i,T_\beta]\right)X_i\\
&=\sum_{i=1}^kg\left(T^\prime_\beta,\sum_{\alpha=1}
^{m-k}\left\{\sum_{j=1}^kA^j_\beta g([X_i,X_j],T_\alpha)+g([X_i,T_\beta],T_\alpha)\right\}T^\prime_\alpha\right)X_i\\
&=\sum_{i=1}^k\left(\sum_{j=1}^kA^j_\beta
\bar{g}([X_i,X_j],T_\beta)+\bar{g}([X_i,T_\beta],T_\beta)\right)X_i
\end{align*}
Thus $\mathcal{P}(\nabla_{T^\prime_\beta}T^\prime_\beta)=0$ if and
only if
 \begin{equation}\label{eq11}
 \sum_{j=1}^kA^j_\beta
\bar{g}([X_i,X_j],T_\beta)=-\bar{g}(\bar{F}_\beta,X_i)
 \end{equation}
for any $i=1,\cdots,k$. The first part is from elementary knowledge
of linear algebra. Since $C_{ij}^\beta=\bar{g}([X_i,X_j],T_\beta)$,
the uniqueness follows from \eqref{eq11} and the assumption that the
matrix $[C_{ij}^\beta]$ is invertible for any $\beta$.
\end{proof}
\begin{proposition}\label{strongbracket} Let $\Sigma$ be a
distribution of $M$. Then $\Sigma$ is strong-bracket generating  if
and only if for any local frame of $TM$, $\{X_1,\cdots,X_k,$
$T_1,\cdots,T_{m-k}\}$, where $\{X_1,\cdots,X_k\}$ is a basis of
$\Gamma(\Sigma)$, the matrix $[C_{ij}^\beta]$ is invertible for
$\beta=1,\cdots,m-k,$ where $C_{ij}^\beta$ are the coefficients
satisfying \eqref{coefficients}.
\end{proposition}
\proof Denote by $\Sigma^\mathcal{\bot}$ be the set of all sections
in $T^\ast M$ annihilating $\Sigma$. By the Cartan formula
$$
d\omega(X,Y)=X(\lambda(Y))-Y(\lambda(X))-\lambda([X,Y]),
$$
it is easy to verify that $\Sigma$ is strong-bracket generating if
and only if
$d\omega:\Gamma(\Sigma)\times\Gamma(\Sigma)\rightarrow\mathbb{R}$ is
nondegenerate for any $\omega\in{\Sigma^\mathcal{\bot}}$, see e.g.
\cite[p.70]{Montgomery}. Now assume $\Sigma$ be strong-bracket
generating. For a given frame $\{X_1,\cdots,X_k,$
$T_1,\cdots,T_{m-k}\}$ of $TM$, for $\beta=1,\cdots,m-k$ we choose
$\omega^\beta\in{\Sigma^\mathcal{\bot}}$ such that
$\omega^\beta(T_\alpha)=\delta_\alpha^\beta$. Then the nondegeneracy
 of $\omega^\beta$ implies the matrix
 $[C_{ij}^\beta]=-[\lambda^\beta([X_i,X_j])]$ is invertible.

Conversely if $0\neq\omega\in{\Sigma^\mathcal{\bot}}$, then we can
choose a frame $\{X_1,\cdots,X_k,$ $T_1,\cdots,T_{m-k}\}$  such that
$\{X_1,\cdots,X_k\}$ is a basis of $\Gamma(\Sigma)$, $\omega(T_1)=1$
and $\omega(T_\beta)=0$ for $\beta\neq1$. The nondegeneracy of
$[C_{ij}^1]$ implies the nondegeneracy of $\omega$ on
$\Gamma(\Sigma)$.
\endproof

 If we write out $\Delta^\mathcal{H}$ in terms of horizontal
vector fields, we see that $\Delta^\mathcal{H}$ is a H\"{o}rmander
operator. In fact, we have
\begin{lemma}\label{laplacian}Let $\{X_i\}_{i=1}^k$ be any
 orthonormal basis of $\Sigma$. Then
$ \Delta^\mathcal{H}=\sum_{i=1}^k(X_i^2-D_{X_i}X_i). $
 \end{lemma}
 Thus by \cite{Hor} $\Delta^\mathcal{H}$ is hypoelliptic if $M$ is connected and $\Sigma$
 satisfies the H\"{o}rmander condition. We will
 use the following technical lemma.
\begin{lemma}\label{technical} Let $M,\Sigma,\Sigma^{\prime},g,g_c$ be as in Theorem
\ref{main1}. For $\epsilon>0$, let $g^\epsilon$ be the Riemannian
metric $ g^\epsilon=g_c\bigoplus \epsilon^2g^\prime $ where
$g^{\prime}:=g|_{\Sigma^{\prime}}$. Denote by $\Delta^\epsilon$ the
Riemannian Laplacian of $g^\epsilon$. Then
$$\lim_{\epsilon\rightarrow+\infty}-\Delta^\epsilon\\=-\Delta^{\mathcal{H}}+\mathrm{H}^\bot
$$
\end{lemma}
\proof Denote by $\nabla^\epsilon$ the Riemannian connection of
$g^\epsilon$. Assume $\{X_1,\cdots,X_k,T_1,\cdots,$ $T_{m-k}\}$ an
orthonormal basis with respect of $g$. Then
$\{X_1,\cdots,X_k,\frac{1}{\epsilon}T_1,\cdots,\frac{1}{\epsilon}T_{m-k}\}$
is an orthonormal basis with respect to $g^\epsilon$. For any smooth
function $f$, by definition and Lemma \ref{laplacian} we get
\begin{align*}
\Delta^\epsilon
f&=\sum_{i=1}^k\left(X_i^2-\nabla^\epsilon_{X_i}X_i\right)f+\frac{1}{\epsilon^2}
\sum_{\beta=1}^{m-k}\left(T_\beta^2-\nabla^{\epsilon}_{T_\beta}T_\beta\right)f\\
&=\Delta^\mathcal{H}f+\sum_{i=1}^kB(X_i,X_i)f+\frac{1}{\epsilon^2}
\sum_{\beta=1}^{m-k}\left(T_\beta^2-\nabla^{\epsilon}_{T_\beta}T_\beta\right)f,
\end{align*}
where
$$
B(X_i,X_i)=\sum_{\beta=1}^{m-k}g^\epsilon\left(\nabla^\epsilon_{X_i}X_i,\frac{1}{\epsilon}T_\beta\right)\frac{1}
{\epsilon}T_\beta=\frac{1}{\epsilon^2}\sum_{\beta=1}^{m-k}g_c(X_i,[T_\beta,X_i]^\mathcal{H})T_\beta
$$
and
\begin{align*}
\nabla^\epsilon_{T_\beta}T_\beta&=(\nabla^\epsilon_{T_\beta}T_\beta)^\mathcal{H}+(\nabla^\epsilon_{T_\beta}T_\beta)^\bot\\
&=\sum_{i=1}^kg^\epsilon(\nabla^\epsilon_{T_\beta}T_\beta,X_i)X_i+\sum_{s=1}^{m-k}g^\epsilon\left(\nabla^\epsilon_{T_\beta}T_\beta,\frac{1}
{\epsilon}T_s\right)\frac{1}{\epsilon}T_s\\
&=
\epsilon^2\sum_{i=1}^kg(T_\beta,[X_i,T_\beta])X_i+\frac{1}{\epsilon^2}\sum_{s=1}^{m-k}g^\epsilon
\left(\nabla^\epsilon_{T_\beta}T_\beta,T_s\right)T_s\\
&=\epsilon^2\mathrm{H}^\bot+\sum_{s=1}^{m-k}g(T_\beta,[T_s,T_\beta])T_s.
\end{align*}
\endproof
\begin{remark} (1),
Some authors called
$\bar{\Delta}^\mathcal{H}:=\Delta^\mathcal{H}-\mathrm{H}^\bot$
sublaplacian, \cite{Fu,Ge}. If $\mathrm{H}^\bot\neq0$,
$\bar{\Delta}^\mathcal{H}$ explicitly depends on $g^\prime$. (2),
the penalty metric $g^\epsilon$ is very useful in sub-Riemannian
geometry. The reason is that when $\Sigma$ satisfies the
H\"{o}rmander condition and $M$ is connected, $(M,d^\epsilon)$
($d^\epsilon$ is the Riemannian distance corresponding to
$g^\epsilon$) converges to $(M,d_{cc})$ in the sense of
Hausdorff-Gromov, e.g. \cite{Mi,Ger,Montgomery}.
\end{remark}
\proof[\textbf{Proof of Theorem \ref{main1}}]Let $g^\epsilon$ as in
Lemma \ref{technical}. Denote by $d\mathrm{(vol)}^\epsilon$ the
Riemannian volume element of $g^\epsilon$. It is easy to show that
\begin{equation}\label{volume}
d\mathrm{(vol)}^\epsilon=\epsilon^{m-k}d\mathrm{vol}.
\end{equation}
Let $e,f$ be any smooth (or Sobolev) functions with compact support.
We abuse the notation to denote by $\nabla^\epsilon f$ the
Riemannian gradient of $f$ with respect to $g^\epsilon$. Noting that
$$
\nabla^\epsilon
f=\nabla^\mathcal{H}f+\frac{1}{\epsilon^2}\sum_{\beta=1}^{m-k}(T_\beta
f)T_\beta
$$
and hence
$$
g^\epsilon(\nabla^\epsilon e,\nabla^\epsilon
f)=g_c(\nabla^\mathcal{H}e,\nabla^\mathcal{H}f)+\frac{1}{\epsilon^2}\sum_{\beta=1}^{m-k}(T_\beta
e)(T_\beta f),
$$
from \eqref{volume} and the Green formula
$$
\int_M(-\Delta^\epsilon
e)fd\mathrm{(vol)}^\epsilon=\int_Mg^\epsilon(\nabla^\epsilon
e,\nabla^\epsilon f)d\mathrm{(vol)}^\epsilon
$$
we derive
\begin{equation}\label{eq1}
\int_M (-\Delta^\epsilon e)fd\mathrm{vol}=\int_M\left\{
g_c(\nabla^\mathcal{H}e,\nabla^\mathcal{H}f)+\frac{1}{\epsilon^2}\sum_{\beta=1}^{m-k}(T_\beta
e)(T_\beta f)\right\}d\mathrm{vol}.
\end{equation}
Taking the limit $\epsilon\rightarrow+\infty$ in \eqref{eq1},  we by
Lemma \ref{technical} induce
\begin{equation}\label{eq2}
    \int_M \left((-\Delta^\mathcal{H}+\mathrm{H}^\bot)e\right)fd\mathrm{vol}=\int_M
g_c(\nabla^\mathcal{H}e,\nabla^\mathcal{H}f)d\mathrm{vol}.
\end{equation}
Putting $f=u\bar{f}(\bar{f}\in{C_0^\infty})$ in \eqref{eq2}, we get
$$
-\int_M \Delta^\mathcal{H}e\bar{f}ud\mathrm{vol}+\int_M
(\mathrm{H}^\bot e)\bar{f}ud\mathrm{vol}=\int_M
g_c(\nabla^\mathcal{H}e,\nabla^\mathcal{H}\bar{f})ud\mathrm{vol}+\int_M
g_c(\nabla^\mathcal{H}e,\nabla^\mathcal{H}u)\bar{f}d\mathrm{vol}.
$$
Thus
\begin{equation}\label{eq3}
-\int_M (\Delta^\mathcal{H}e)\bar{f}d\mu=\int_M
g_c(\nabla^\mathcal{H}e,\nabla^\mathcal{H}\bar{f})d\mu,
\end{equation}
holds if and only if
$$
\int_M (\mathrm{H}^\bot e)\bar{f}ud\mathrm{vol}=\int_M
g_c(\nabla^\mathcal{H}e,\nabla^\mathcal{H}u)\bar{f}d\mathrm{vol}.
$$
By the arbitrariness of $\bar{f}$ in the last equation, we deduce
that \eqref{eq3} holds for any $f,\bar{f}\in{C_0^\infty(M)}$ if and
only
$$
(\mathrm{H}^\bot e)u=g_c(\nabla^\mathcal{H}e,\nabla^\mathcal{H}u)
$$
for any $e\in{C_0^\infty(M)}$, that is,
$$
\mathrm{H}^\bot e=g_c(\nabla^\mathcal{H}e,\nabla^\mathcal{H}(\ln
u)).
$$
Since $\mathrm{H}^\bot$ is a horizontal vector field and hence
$\mathrm{H}^\bot e=g_c(\nabla^\mathcal{H}e,\mathrm{H}^\bot)$,
Theorem \ref{main1} follows.
\endproof
\begin{corollary}\label{symmetric}
Let $(M,\Sigma,g_c)$ be a sub-Riemannian manifold. Then there exists
a complement $\Sigma^\prime$ of $\Sigma$ such that we can extend
$g_c$ to some Riemannian metric $g$ and $\Delta^\mathcal{H}$ is a
symmetric operator on $C_0^\infty(M)$: for any
$e,f\in{C_0^\infty(M)}$
$$
\int_M (-\Delta^\mathcal{H}e) fd\mathrm{vol}=\int_M
g_c(\nabla^\mathcal{H}e,\nabla^\mathcal{H}f)d\mathrm{vol}=\int_M e
(-\Delta^\mathcal{H}f)d\mathrm{vol}
$$
where $d\mathrm{vol}$ is the Riemannian measure of $g$.
\end{corollary}
\begin{definition}\label{defn}Let $(M,\Sigma,g_c)$ be a sub-Riemannian manifold.
Fix a complement of $\Sigma$ such that \eqref{curvature} holds for
some extension $g$ of $g_c$ and for some orthonormal basis of $g$.
We define $\Delta^\mathcal{H}$ (with respect to the splitting
$TM=\Sigma\bigoplus\Sigma^\prime$) as \textbf{a sublaplacian of}
$\mathbf{(M,\Sigma,g_c)}$. When such complement is unique (see
Theorem \ref{main2}), we call $\Delta^\mathcal{H}$ \textbf{the
sublaplacian of} $\mathbf{(M,\Sigma,g_c)}$.
\end{definition}
One of the reasons we define $\Delta^{\mathcal{H}}$ (not
$\bar{\Delta}^\mathcal{H}$) as a (the) sublaplacian is that
$\Delta^{\mathcal{H}}$ is compatible with several notions such as
the horizontal Hessian and weakly convex functions on sub-Riemannian
manifolds, see \cite{T}. As already pointed out in the introduction,
\textbf{the laplacian} in Definition \ref{defn} is defined for few
cases. Theorem \ref{main2} and Proposition \ref{strongbracket} tell
us that the sublaplacian is well defined for fat sub-Riemannian
manifolds. This make the case more interesting because fat
sub-Riemannian manifolds are proven to admit no singular
sub-Riemannian geodesics.
\begin{example}[Carnot groups, \cite{FS,Pansu}]A Carnot group (or a stratified group)
$G$ is a connected, simply connected Lie group whose Lie algebra
$\mathcal{G}$ admits the grading
$\mathcal{G}=V_{1}\bigoplus\cdots\bigoplus V_{l}$, with
$[V_{1},V_{i}]=V_{i+1}$, for any $1\leq i\leq l-1$ and
$[V_{1},V_{l}]={0}$ (the integer $l$ is called the step of $G$). Let
$\{e_{1},\cdots,e_{m}\}$ be a basis of $\mathcal{G}$ with
$m=\sum_{i=1}^{l}\dim(V_{i})$. Let $X_{i}(x)=(L_{x})_{*}e_{i}$ for
$i=1,\cdots,k:=\dim(V_{1})$ where $(L_{x})_{*}$ is the differential
of the left translation $L_{x}(x^{\prime})=xx^{\prime}$ and let
$T_{\alpha}(x)=(L_{x})_{*}e_{i+k}$ for $\alpha=1,\cdots,m-k$. We
call the system of left-invariant vector fields
$\Sigma:=V_{1}=\textrm{span}\{X_{1},\cdots,X_{k}\}$ the horizontal
bundle of $G$. If we equip $\Sigma$ an inner product $g_{c}$ such
that $\{X_{1},\cdots,X_{k}\}$ is an orthonormal basis of $\Sigma$,
$(G,\Sigma,g_{c})$ is a sub-Riemannian manifold satisfying the
H\"{o}rmander condition which is guaranteed by the grading of its
Lie algebra. The role played by Carnot groups in sub-Riemannian
geometry is similar that by Euclidean spaces in Riemannian geometry,
\cite{Mi}. Thus sub-Riemannian manifolds are also called Carnot
manifolds. Fix a Carnot group $G$. Because of the grading condition
of its Lie algebra, by choosing the natural splitting of $TG$ and a
system of left-invariant vector fields $\{X_1,\cdots,X_k\}$ as an
orthonormal basis, we easily deduce that the horizontal connection
$D$ has the following simple form
$$
D_XY=\sum_{i=1}^kX(Y^i)X_i,\quad\mathrm{ for\; any \quad
}X,Y=\sum_{i=1}^kY^iX_i\in{\Gamma(\Sigma)}
$$
and hence
$$
\Delta^\mathcal{H}=\sum_{i=1}^kX_{i}^2.
$$
$\Delta^\mathcal{H}$ coincides with the sublaplacian of Carnot
groups studied in the literature, see \cite{RS,Je,BU} and references
therein. It is clear  that $\Delta^\mathcal{H}$ is a symmetric
operator on $C_0^\infty(G)\hookrightarrow L^2(G,g)$, where $g$ is an
extension of $g_c$ such that $\{X_1,\cdots,X_k,T_1,\cdots,T_{m-k}\}$
is  orthonormal.
\end{example}
\begin{example}[contact Riemannian
manifolds,\cite{We,Bl}] Let $M$ be a real $2n + 1$-dimensional
smooth manifold. An almost contact Riemannian structure
$(\varphi,\xi,\eta,g)$ on $M$ consists of a $(1, 1)$- tensor field
$\varphi$, a vector field $\xi$, a 1-form $\eta$, and a Riemannian
metric $g$ such that
$$
\varphi^2=-I+\eta\otimes\xi,\quad\eta\circ\varphi=0,\varphi\xi=0,\quad
g(\varphi X,\varphi Y)=g(X,Y)-\eta(X)\eta(Y),
$$
for any $X,Y\in{\Gamma(M)}$. It is a contact Riemannian structure if
it satisfies $\Omega=d\eta$ (the contact condition) where $\Omega(X,
Y) = g(X,\varphi Y)$. Set $\Sigma=\ker(\eta)$. Then
$(M,\Sigma,g|_\Sigma)$ is a sub-Riemannian structure and $\Sigma$
has a natural complement. Let $\nabla$ be the Levi-Civita connection
of $g$. Then the horizontal connection $D$ is
$$
D_XY=\nabla_XY-\eta(\nabla_XY)\xi, \quad X,Y\in{\Gamma(\Sigma)}.
$$
Note that the Levi form
$$
L_\eta(X,Y)=-d\eta(X,\varphi Y)=\eta([X,\varphi Y])\quad
X,Y\in{\Gamma(\Sigma)}
$$ is nondegenerate. This in particular implies that $\Sigma$
satisfies the H\"{o}rmander condition. The (generalized)
Tanaka-Webster connection (\cite{Ta}) on $(M,\varphi,\xi,\eta,g)$ is
$$
\nabla^\ast_XY=\nabla_XY+\eta(X)\varphi(Y)-
\eta(Y)\nabla_X\xi+[(\nabla_X\eta)Y]\xi.
$$
Denote by $\mathcal {D}$ the complexification of $\Sigma$, i.e.,
$\mathcal{D}=\mathcal{D}^\prime\oplus\mathcal{\overline{D}}^\prime$
where $\mathcal{D}^\prime=\{X-i\varphi X,X\in{\Gamma(\Sigma)}\}$ and
$\mathcal{\overline{D}}^\prime$ is the conjugate of
$\mathcal{D}^\prime$.  Set $h:=\frac{1}{2}\mathcal{L}_\xi\varphi$.
Then the pair $(M,\mathcal{D})$ is a (strongly pseudo-convex) CR
manifold, (i.e., $[\mathcal{D}^\prime,\mathcal{D}^\prime]\subset
\mathcal{D}^\prime$ and the Levi form $L_\eta$ is positive
definite), if and only if the contact Riemannian manifold
$(M,\varphi,\xi,\eta,g)$ satisfies
$$
(\nabla_X\varphi)Y=g(X+hX,Y)\xi-\eta(Y)(X+hX).
$$
Let $(M,\mathcal{D})$ be a strongly pseudo-convex CR manifold.
Denote by $D^\ast$ the restriction on $\Sigma$ of the $\nabla^\ast$,
and extend $D^\ast$ to the complexified bundle $\mathcal{D}$. Note
that $D^\ast$ is just the Webster connection, \cite{We}. Since for
$X,Y\in{\Gamma(\Sigma)}$, $\eta(X)=\eta(Y)=0$,
$(\nabla_X\eta)Y=X(\eta(Y))-\eta(\nabla_XY)$, we have $D^\ast=D$.
Here $D$  is also extended to $\mathcal{D}$. Let
$\{X_\alpha\}_{\alpha=1}^n$ be an orthonormal complex basis (with
respect to the extended metric $g$) of $\mathcal{D}^\prime$. Then
$\{X_{\overline{\alpha}}=\varphi(X_\alpha)\}_{\alpha=1}^n$ is an
orthonormal complex basis of $\mathcal{\overline{D}}^\prime$. For a
smooth function $f$ on $M$, the sublaplacian studied by \cite{Gr,LL}
is
$$
\overline{\Delta}f=\sum_{\alpha=1}^nf_{\alpha\overline{\alpha}}+f_{\overline{\alpha}\alpha}
$$
where
$$
f_\alpha=X_\alpha
f,\:f_{\overline{\alpha}}=X_{\overline{\alpha}}f,\quad
f_{\alpha\beta}=X_\beta
f_\alpha-\sum_{\gamma=1}^n\Gamma_{\alpha\beta}^\gamma f_\gamma,\quad
f_{\overline{\alpha}\beta}=X_\beta
f_{\overline{\alpha}}-\sum_{\gamma=1}^n\Gamma_{\overline{\alpha}\beta}^{\overline{\gamma}}f_{\overline{\gamma}}
$$
and
$\Gamma_{\alpha\beta}^\gamma=g(D^\ast_{X_\alpha}X_\beta,X_\gamma),\;\Gamma_{\overline{\alpha}\beta}^{\overline{\gamma}}
=g(D^\ast_{X_{\overline{\alpha}}}X_\beta,X_{\overline{\gamma}})$.
Set
$Y_\alpha=\frac{1}{\sqrt{2}}(X_\alpha+X_{\overline{\alpha}}),\overline{Y}_\alpha=\frac{i}{\sqrt{2}}(X_\alpha-X_{\overline{\alpha}})$.
Then $\{Y_1,\cdots,Y_n,\overline{Y}_1,\cdots,\overline{Y}_n\}$ is an
orthonormal basis of $\Sigma$. Now by direct computation we get from
Lemma \ref{laplacian}
$$
\overline{\Delta}=\sum_{\alpha=1}^n(Y^2_\alpha+\overline{Y}^2_\alpha-
D_{Y_\alpha}Y_\alpha-D_{\overline{Y}_\alpha}\overline{Y}_\alpha)=\Delta^\mathcal{H},
$$
since $D^\ast=D$. Thus for strongly pseudo-convex pseudo-Hermitian
manifolds our definition for sublaplacians coincides with the
canonical one. Because $\nabla_\xi \xi=0$ and the Levi form $L_\eta$
is nondegenerate, by Corollary \ref{symmetric}, $\Delta^\mathcal{H}$
is a symmetric operator on $C_0^\infty(M)\hookrightarrow L^2(G,g)$.
\end{example}
\begin{example}[Riemannian submersions with minimal fibres, \cite{FLP}] Let $(M,g)$ and $(B,g^\prime$
 be smooth Riemannian manifolds. A smooth map $\pi:M\rightarrow B$ is a
 submersion if $\pi_\ast:T_qM\rightarrow T_{\pi(q)}B$ is a
 surjective linear map for each $q\in{M}$. The vertical space at $q$
 is the tangent space of the fibre $\pi^{-1}(\pi(q))$:
 $V_q=\ker(\pi_\ast)$. The collection of vertical spaces is the
 vertical distribution $V\subset TM$. Let $\Sigma$ be the orthogonal
 complement of $V$. $M$ with the structure $(\Sigma,g_c=g|_\Sigma)$
 is a sub-Riemannian manifold. If $\pi_\ast:\Sigma_q\rightarrow
 T_{\pi(q)}$ is linear isometry for any $q\in{M}$, $\pi$ is called a
 Riemannian submersion. The Riemannian submersion $\pi$ is a
 harmonic map between $(M,g)$ and $(B,g^\prime)$ if and only if each
 fibre of $\pi$ is a minimal surface, e.g. \cite{FLP}. If $\pi$
 is a Riemannian submersion with minimal fibres, by Theorem
 \ref{main1} we can define a sublaplacian $\Delta^\mathcal{H}$ on
 $(M,\Sigma,g_c)$ such that $\Delta^\mathcal{H}$ is a symmetric
operator on $C^\infty_0(M)\hookrightarrow L^2(M,g)$.
\end{example}
The above examples show that our notion of sublaplacians covers the
canonical ones in the literature.
\begin{lemma}[divergence theorem]
Let  $(M,\Sigma,g_c)$ be a sub-Riemannian manifold. Let $g$ be the
orthogonal extension of $g_c$ as in Theorem \ref{main2}. Then for
any horizontal vector field $X\in{\Gamma(\Sigma)}$
\begin{equation}\label{divergence}
   \mathrm{div}^\mathcal{H}X=\mathrm{ div}X
\end{equation}
where $\mathrm{div}$ is the Riemannian divergence of $g$. Thus if
$M$ is moreover compact with boundary (possibly empty), we have for
any horizontal vector field $X$
\begin{equation}\label{divergencetheorem}
\int_M \mathrm{div}^\mathcal{H}Xd\mathrm{vol}=\int_{\partial
M}g(X,\nu)ds
\end{equation}
where $d\mathrm{vol}$ is the Riemannian measure of $g$, $\nu$ is the
normal vector field of the boundary $\partial M$, and $ds$ is the
area measure on $\partial M$ induced by $g$.
\end{lemma}
\proof Choose $\{X_1,\cdots,X_k,T^\prime_1,\cdots,T^\prime_{m-k}\}$
as an orthonormal basis of $g$, such that \eqref{curvature} holds.
Since the (horizontal) divergence is independent of the choice of
orthonormal bases,
\begin{align*}
\mathrm{div}X&=\sum_{i=1}^kg(\nabla_{X_i}X,X_i)+\sum_{\alpha=1}^{m-k}g(\nabla_{T^{\prime}_\alpha}X,T^{\prime}_\alpha)\\
&=\sum_{i=1}^kg_c(D_{X_i}X,X_i)+\sum_{\alpha=1}^{m-k}\left\{T^\prime_\alpha
g(X,T^\prime_\alpha)-g(X,\nabla_{T^\prime_\alpha}T^\prime_\alpha)\right\}\\
&=\mathrm{div}^\mathcal{H}
\end{align*}
where we used \eqref{curvature} and the assumption that $X$ is
horizontal. \eqref{divergencetheorem} is from \eqref{divergence} and
the classical divergence theorem.
\endproof
\begin{theorem}Let $\Delta^\mathcal{H}$ be a (the) sublaplacian of $(M,\Sigma,g_c)$ in
the sense of Definition \ref{defn}.  If $\Sigma$ satisfies the
H\"{o}rmander condition and $M$ is a closed, connected manifold,
then any horizontal-harmonic function $u$, i.e. $u$ satisfies
$$
\Delta^\mathcal{H}u=0,
$$
is constant.
\end{theorem}
\proof Note that
\begin{equation}\label{eq12}
\Delta^\mathcal{H}u^2=2\mathrm{div}^\mathcal{H}(u\nabla^\mathcal{H}u)
=2g_c(\nabla^\mathcal{H}u,\nabla^\mathcal{H}u)+2u\Delta^\mathcal{H}u
\end{equation}
If $u$ is horizontal-harmonic, by \eqref{divergence} and
\eqref{eq12} we get
$$
\mathrm{div}(u\nabla^\mathcal{H}u)=g_c(\nabla^\mathcal{H}u,\nabla^\mathcal{H}u)
$$
where $\mathrm{div}$ is the Riemannian divergence of some extension
$g$ of $g_c$ as in Theorem \ref{main2}. Integrating the last
formula, by the green formula in the Riemannian case we induce
$$
\int_Mg_c(\nabla^\mathcal{H}u,\nabla^\mathcal{H}u)d\mathrm{vol}=0
$$
since by assumption $M$ is closed. Thus $\nabla^\mathcal{H}u=0$,
that is, $u$ is constant along horizontal curves. The statement
follows because $\Sigma$ satisfies the H\"{o}rmander condition and
$M$ is connected, by the Chow theorem \cite{Chow} any two points can
be connected by a piecewisely smooth horizontal curve.
\endproof
\section{Eigenvalues of sublaplacians of compact sub-Riemannian manifolds}
In this section we always assume $(M,\Sigma,g_c)$ is a
\textbf{compact and regular} sub-Riemannian manifold with smooth
(possibly empty) boundary. Let $\Delta^\mathcal{H}$ be a (the)
sublaplacian of $(M,\Sigma,g_c)$ in the sense of Definition
\ref{defn} and $g$ be an orthogonal extension of $g_c$ with respect
to the given decomposition as in Definition \ref{defn}. The goal of
this section is to study the eigenvalue problem of
$\Delta^\mathcal{H}$. First we give the definition of horizontal
Sobolev functions on $(M,\Sigma,g_c)$.
\begin{definition}A function $f$ in $L^2(M)$ is called a horizontal
Sobolev function if there exists a horizontal vector filed $Y$
belonging to $L^2(M)$ such that the following
$$
\int_Mg_c(Y,X)d\mathrm{vol}=-\int_Mf\mathrm{div}^\mathcal{H}Xd\mathrm{vol}
$$
holds for any horizontal vector field $X$ with compact support on
$M$. $Y$ denoted by $\nabla^\mathcal{H}f$ is called the weakly
horizontal derivative of $f$. The set of all horizontal Sobolev
functions is denoted by $W^{1,2}_{\mathcal{L}}(M)$.
\end{definition}
Here we call a horizontal vector field is in $L^2(M)$ if its
coefficients are in $L^2(M)$. From \eqref{divergencetheorem} the
above definition is well-defined. We denote by
$H^{1,2}(M)(H_0^{1,2}(M))$ the completed space of
$C^\infty(M)(C_0^\infty(M))$ functions with respect to the norm
$$
||f||_{H^{1,2}(M)}=\left(\int_M|f|^2+g_c(\nabla^\mathcal{H}f,\nabla^\mathcal{H}f)d\mathrm{vol}\right)^\frac{1}{2}.
$$
\begin{lemma}\label{sobolev}\quad\\
(1),$H^{1,2}(M)=W^{1,2}_{\mathcal{L}}(M);$\\
(2), the embedding $W^{1,2}_{\mathcal{L}}(M)\hookrightarrow L^2(M)$ is compact;\\
(3), If $f\in{W^{1,2}_{\mathcal{L}}(M)}$ and
$\Delta^\mathcal{H}f=\lambda f$ for some $\lambda\in{\mathbb{R}}$,
then $f$ must be smooth.
\end{lemma}
\proof Since $M$ is compact, by choosing a smooth partition of unity
subordinate to a \textbf{finite} cover of $M$, the first two
statements are reduced to a local chart case. Let $\phi:U\subset
M\rightarrow V=\phi(U)\subset\mathbb{R}^m$ be a coordinate chart.
Since $\Sigma$ is regular, $\phi_\ast(\Sigma|_U)$ is also a regular
distribution on $V$. Then $(V,\phi_\ast(\Sigma|_U),g^\prime)$ is a
regular sub-Riemannian manifold, where $g^\prime$ is the standard
Euclidean metric. Now any function in $H^{1,2}(U)$ is pushed forward
by $\phi$ to a horizontal weighted Sobolev function in the sense of
\cite{FSSC1}. Now the first two statements follow from the
corresponding results proven in \cite{FSSC1,FLW2}. The third is
standard since $\Delta^\mathcal{H}$ is a hypoelliptic operator.
\endproof
\begin{theorem}\label{spectrsl}Let $M$ be without boundary. Consider the following
closed eigenvalue problem
\begin{equation}\label{eigen}
-\Delta^\mathcal{H}f=\lambda f.
\end{equation}
That is, we are looking for all numbers $\lambda$ for which there
exists a nontrivial smooth solution satisfying \eqref{eigen}. Then
\begin{enumerate}
  \item The set of eigenvalues consists of an infinite sequence $0\leq\lambda_1<\lambda_2<\lambda_3\cdots\uparrow+\infty$
  \item Each eigenvalue $\lambda_i$ has finite multiplicity and the
  eigenspaces corresponding to different eigenvalues are $L^2(M)-$
  orthogonal.
  \item The direct sum of the eigenspaces $E(\lambda_i)$,
  $i=1,\cdots$ is dense in $L^2(M)$.
  \item Let $\Delta^\epsilon$ be as in Lemma \ref{technical}. For
  each $\epsilon$, denote by $\lambda_i(\epsilon)$ be the $i$-th (counting the
  multiplicity)   eigenvalue of the eigenvalue problem
  $$
  \Delta^\epsilon f=\lambda(\epsilon)f.
  $$
  Then
  $$
  \lim_{\epsilon\rightarrow+\infty}\lambda_i(\epsilon)=\lambda_i
  $$
  where $\lambda_i$ is renumbered counting the multiplicity.
\end{enumerate}
\end{theorem}
\proof By Corollary \ref{symmetric}, $-\Delta^\mathcal{H}$ is a
positive and symmetric operator on $C^\infty(M)$ which is dense in
$W_\mathcal{L}^{1,2}(M)$. Thus by the first statement of Lemma
\ref{sobolev}, $-\Delta^\mathcal{H}$ can be extended to a closed,
positive self-adjoint operator on $W_\mathcal{L}^{1,2}(M)$, which
implies that the spectrum of $-\Delta^\mathcal{H}$ is contained in
$\mathbb{R}_+$. It follows from the compactness of the embedding
$W^{1,2}_{\mathcal{L}}(M)\hookrightarrow L^2(M)$ that the resolvent
$(-\Delta^\mathcal{H}-\lambda)^{-1}$ is a compact operator in
$L^2(M)$. The first three statement follows from the classical
results on the spectral theory of compact operators and from the
third claim of Lemma \ref{sobolev}, see e.g. \cite{DS}.

Fukaya in \cite{Fu} proved the fourth statement for
$\Delta^\mathcal{H}-\mathrm{H}^\bot$. Since by our choice of the
orthogonal extension  $\mathrm{H}=0$, the statement follows.
\endproof
\begin{remark}(1),The first three claims can also be proven by a
variational argument minimizing the Rayleigh quotient
$$
\frac{\int_M
g_c(\nabla^\mathcal{H}f,\nabla^\mathcal{H}f)d\mathrm{vol}}{\int_Mf^2d\mathrm{vol}}.
$$
(2),For complete sub-Riemannian manifolds, following the lines of
\cite{Str} it is possible to develop a theory of heat semi-group of
$\Delta^\mathcal{H}$.
\end{remark}

\end{document}